\documentclass[14pt]{article}
\usepackage{mathrsfs}
\usepackage{amsthm}
\usepackage{amssymb}
\usepackage{amsmath}
\usepackage{graphicx}
\usepackage{color}
\usepackage{amsfonts}
\usepackage{float}
\usepackage[text={140mm,210mm},left=45mm,vmarginratio=1:1]{geometry}
\newtheorem{theorem}{Theorem}[section]

\newtheorem{lemma}[theorem]{Lemma}

\numberwithin{equation}{section}
\normalsize

\begin{document}
\title{\textbf{The contact process with semi-infected state on the complete graph}}

\author{Xiaofeng Xue \thanks{\textbf{E-mail}: xfxue@bjtu.edu.cn \textbf{Address}: School of Science, Beijing Jiaotong University, Beijing 100044, China.}\\ Beijing Jiaotong University}

\date{}
\maketitle

\noindent {\bf Abstract:}In this paper we are concerned with the
contact process with semi-infected state on the complete graph $C_n$
with $n$ vertices. In our model, each vertex is in one of three
states that `healthy', `semi-infected' or `wholly-infected'. Only
wholly-infected vertices can infect others. A healthy vertex becomes
semi-infected when being infected while a semi-infected vertex
becomes wholly-infected when being further infected. Each (semi- and
wholly-) infected vertex becomes healthy at constant rate. Our main
result shows the phase transition for the time wholly-infected
vertices wait for to die out. Conditioned on all the vertices are
wholly-infected when $t=0$, we show that wholly-infected vertices
survive for $\exp\{O(n)\}$ units of time when the infection rate
$\lambda>4$ while die out in $O(\log n)$ units of time when
$\lambda<4$.

\quad


\section{Introduction}\label{section one}
In this paper we are concerned with the contact process with
semi-infected state on the complete graph. A complete graph is a
finite graph such that for any two vertices there is an edge
connecting them. For later use, for integer $n\geq 1$, we denote by
$C_n$ the complete graph with $n$ vertices and denote by
$\{1,2,\ldots,n\}$ the vertices set of $C_n$.

The contact process with semi-infected state on $C_n$ is a
continuous-time Markov process with state space $\{0,1,2\}^{C_n}$,
i. e. at each vertex there is a spin taking value from $\{0,1,2\}$.
For any configuration $\eta\in \{0,1,2\}^{C_n}$ and $1\leq i\leq n$,
we denote by $\eta(i)$ the value of the spin at the vertex $i$. For
any $t\geq 0$, we denote by $\eta_t$ the configuration of our
process at moment $t$. For $\eta\in \{0,1,2\}^{C_n}$, $1\leq i\leq
n$ and $0\leq l\leq 2$, we define $\eta^{i,l}\in \{0,1,2\}^{C_n}$ as
follows.
\[
\eta^{i,l}(j)=
\begin{cases}
\eta(j), &\text{~if~}j\neq i,\\
l,& \text{~if~}j=i.
\end{cases}
\]
The generator $\Omega$ of $\{\eta_t\}_{t\geq 0}$ has the form
\begin{equation}\label{equ 1.1 generator}
\Omega f(\eta)=\sum_{1\leq i\leq n}\sum_{l=0,1,2}H(\eta,i,l)\big[f(\eta^{i,l})-f(\eta)\big]
\end{equation}
for any continuous function $f$ on $\{0,1,2\}^{C_n}$. That is to
say, at each moment that the configuration of the process jumps,
only one spin changes value. Conditioned on the current
configuration $\eta$, the process jumps to $\eta^{i,l}$ at rate
$H(\eta,i,l)$. $H(\eta,i,l)$ is defined as
\[
H(\eta,i,l)=
\begin{cases}
1 & \text{~if~}l=0\text{~and~}\eta(i)\neq 0,\\
\frac{\lambda}{n}|\{j:~\eta(j)=2\}| & \text{~if~}l=2\text{~and~}\eta(i)=1,\\
\frac{\lambda}{n}|\{j:~\eta(j)=2\}| & \text{~if~}l=1\text{~and~}\eta(i)=0,\\
0 &\text{~else},
\end{cases}
\]
where $\lambda$ is a positive constant called the infection rate and $|A|$ is the cardinality of the set $A$.

Intuitively, the process describes the spread of an epidemic, where
each individual is in one of the three states that healthy,
semi-infected or wholly-infected. For details, vertices in state $0$
are healthy while vertices in state $1$ are semi-infected and
vertices in state $2$ are wholly-infected. Each infected vertex, no
matter semi- or wholly-, waits for an exponential time with rate $1$
to become healthy. Wholly-infected vertices have the ability to
infect others. A healthy vertex is infected at rate proportional to
the number of wholly-infected vertices and becomes a semi-infected
vertex when being infected. A semi-infected vertex is further
infected at rate proportional to the number of wholly-infected
vertices and becomes a wholly-infected vertex when being infected.

According to the spatial homogeneity of the process
$\{\eta_t\}_{t\geq 0}$, we only care about the numbers of vertices
in state $1$ and $2$. Hence we define
\[
B_t=|\{i:~\eta_t(i)=2\}| \text{~and~}G_t=|\{i:~\eta_t(i)=1\}|
\]
for any $t\geq 0$. According to the generator of $\{\eta_t\}_{t\geq
0}$ given in Equation \eqref{equ 1.1 generator} and the definition
of $H(\eta,i,l)$, the transition rates function of $\{(B_t,
G_t)\}_{t\geq 0}$ is given by
\[
(B,G)\text{~jumps to~}
\begin{cases}
(B-1,G) & \text{~at rate~}B,\\
(B,G-1) & \text{~at rate~}G,\\
(B+1,G-1) & \text{~at rate~}\frac{\lambda}{n}BG,\\
(B,G+1) &\text{~at rate~}\frac{\lambda}{n}B(n-B-G).
\end{cases}
\]
The process $\{\eta_t\}_{t\geq 0}$ is an extension of the classic
contact process introduced in \cite{Har1974} by Harris, where
vertices are distinguished as healthy ones and infected ones such
that a healthy one is infected at rate proportional to the number of
infected neighbors while an infected one becomes healthy at rate
one. For a detailed survey of the study of the classic contact
process, see Chapter 6 of \cite{Lig1985} and Part one of
\cite{Lig1999}.

In \cite{Pet2011}, Peterson studies the contact process on the
complete graph with vertex-dependent infection rates, containing the
classic contact process as a case. It is shown in \cite{Pet2011}
that there is a critical value $\lambda_c$ of the infection rate
$\lambda$ such that when $\lambda<\lambda_c$ the process dies out
before $O(\log n)$ units of time with high probability while when
$\lambda>\lambda_c$ the process survives for $\exp\{O(n)\}$ units of
time with high probability. We are inspired by \cite{Pet2011} a lot.
Our main result shows that similar phase transition with that in
\cite{Pet2011} occurs for the contact process with semi-infected
state. The precise value of the critical infection rate is also
given. For mathematical details, see Section \ref{section two}.

\section{Main result}\label{section two}
In this section we give the main result of this paper. We care about
the first moment that there is no vertex in state $2$, hence we
define
\[
\tau=\inf\{t:B_t=0\}.
\]
Since vertices in state $1$ can not infect others, after the moment
$\tau$, the epidemic dies out in $\log n$ units of time with high
probability, which depends on the fact that the maximum of $n$
independent exponential times with rate one is with order $\log n$
with high probability. As a result, whether the epidemic survives
for a long time, for instance $\exp\{O(n)\}$ units of time, depends
on $\tau$.

For each $n\geq 1$ and any $\lambda>0$, we denote by $P_{\lambda,n}$
the contact process with semi-infected state on $C_n$ with infection
rate $\lambda$. We write $P_{\lambda,n}$ as $P_{\lambda,n}^{B,G}$
when $(B_0,G_0)=(B,G)$. For given $b,g\in [0,1]$, we write
$(B,G)=(nb,ng)$ when $(B,G)=(\lfloor nb\rfloor,\lfloor ng\rfloor)$
for simplicity. We have the following main result which shows phase
transition for $\tau$.

\begin{theorem}\label{theorem 2.1 main}
When $\lambda>4$, there exists a constant $C=C(\lambda)>0$ such that
\begin{equation}\label{equ 2.1 main}
\lim_{n\rightarrow+\infty}P_{\lambda,n}^{n,0}(\tau>e^{Cn})=1
\end{equation}
while when $\lambda<4$,
\begin{equation}\label{equ 2.2 main}
\lim_{n\rightarrow+\infty}P_{\lambda,n}^{n,0}(\tau<(1+\theta)\log
n)=1.
\end{equation}
for any $\theta>0$.
\end{theorem}

According to Theorem \ref{theorem 2.1 main}, conditioned on all the
vertices are wholly-infected at $t=0$, the wholly-infected vertices
survive for $\exp\{O(n)\}$ units of time when $\lambda>4$ while die
out in $O(\log n)$ units of time when $\lambda<4$.

We are inspired by former references about the contact process on
finite sets to prove Theorem \ref{theorem 2.1 main}. According to
the main results in \cite{Pet2011}, the classic contact process on
the complete graph $C_n$ survives for $\exp\{O(n)\}$ units of time
when $\lambda>1$ while dies out in $O(\log n)$ units of time when
$\lambda<1$. Foxall, Edwards and van den Driessche introduce the
contact process on the complete graph incorporating monogamous
dynamic partnerships in \cite{Fox2016}, where similar phase
transition with that in \cite{Pet2011} is shown and the precise
value of the critical infection rate is given. Durrett and coworkers
study the contact process on the finite lattice $[-N,N]^d$ in
\cite{Dur1988, Dur1988b,Dur1989}. One of their main results is that
the process survives for $\exp\{O(N^d)\}$ units of time when
$\lambda>\lambda_c(d)$ while dies out in $O(\log N)$ units of time
when $\lambda<\lambda_c(d)$, where $\lambda_c(d)$ is the minimum of
the infection rates of the contact process on $\mathbb{Z}^d$ with
which the process survives forever with positive probability.

The proof of Theorem \ref{theorem 2.1 main} is divided into three
sections. In Section \ref{section three}, we introduce an
two-dimensional ODE
\[
\begin{cases}
\frac{d}{dt}b_t=F_1(b_t,g_t)\\
\frac{d}{dt}g_t=F_2(b_t,g_t)
\end{cases}
\]
and show that the solution $(b_t,g_t)$ to this ODE is the mean field
limit of $(\frac{B_t}{n},\frac{G_t}{n})$ for $t$ in any compact area
as $n$ grows to infinity. In Section \ref{section four}, we give the
proof of Equation \eqref{equ 2.1 main}. The proof relies heavily on
the fact that when $\lambda>4$ there exists $(b,g)$ such that
$F_1(b,g)>0$ and $F_2(b,g)>0$. In Section \ref{section five}, we
give the proof of Equation \eqref{equ 2.2 main}. The proof relies
heavily on the fact that $(0,0)$ is the unique equilibrium state of
the ODE when $\lambda<4$ .

For later use, at the end of this section we give a lemma which
shows that the contact process with semi-infected state is monotonic
under a specific partial order on $\mathbb{Z}^2$. For $(B_1, G_1),
(B_2,G_2)\in \mathbb{Z}^2$, we write $(B_1,G_1)\succeq (B_2,G_2)$
when and only when $B_1\geq B_2$ and $B_1+G_1\geq B_2+G_2$. It is
easy to check that $\succeq$ is a partial order on $\mathbb{Z}^2$.
We write $(B_2,G_2)\preceq (B_1, G_1)$ when $(B_1,G_1)\succeq
(B_2,G_2)$. We have the following lemma.
\begin{lemma}\label{lemma 2.2} For j=1,2,
if $\{(B_t^j,G_t^j)\}_{t\geq 0}$ is the contact process with
semi-infected state with initial state $(B_0^j,G_0^j)$ and
$(B_0^1,G_0^1)\succeq (B_0^2,G_0^2)$, then $(B_t^1,G_t^1)\succeq
(B_t^2,G_t^2)$ for any $t\geq 0$ in the sense of coupling.
\end{lemma}
\proof

Conditioned on $(B^1,G^1)\succeq (B^2,G^2)$, we couple the two
processes as follows.
\begin{align*}
&\Big((B^1,G^1),(B^2,G^2)\Big)\rightarrow\\
&\begin{cases} \Big((B^1-1,G^1),(B^2-1,G^2)\Big) & \text{~at rate~}B^2,\\
\Big((B^1-1,G^1),(B^2,G^2)\Big) & \text{~at rate~}B^1-B^2,\\
\Big((B^1,G^1-1),(B^2,G^2-1)\Big) & \text{~at rate~}\min\{G^1,G^2\},\\
\Big((B^1,G^1-1),(B^2,G^2)\Big) & \text{~at
rate~}G^1-\min\{G^1,G^2\},\\
\Big((B^1,G^1),(B^2,G^2-1)\Big) & \text{~at
rate~}G^2-\min\{G^1,G^2\},\\
\Big((B^1+1,G^1-1),(B^2+1,G^2-1)\Big) & \text{~at
rate~}\frac{\lambda}{n}\min\{B^1G^1,B^2G^2\},\\
\Big((B^1+1,G^1-1),(B^2,G^2)\Big) & \text{~at
rate~}\frac{\lambda}{n}(B^1G^1-\min\{B^1G^1,B^2G^2\}),\\
\Big((B^1,G^1),(B^2+1,G^2-1)\Big) & \text{~at
rate~}\frac{\lambda}{n}(B^2G^2-\min\{B^1G^1,B^2G^2\}),\\
\Big((B^1,G^1+1),(B^2,G^2+1)\Big) & \text{~at rate~}a,\\
\Big((B^1,G^1+1),(B^2,G^2)\Big) & \text{~at rate~}\frac{\lambda}{n}B^1(n-B^1-G^1)-a,\\
\Big((B^1,G^1),(B^2,G^2+1)\Big) & \text{~at
rate~}\frac{\lambda}{n}B^2(n-B^2-G^2)-a,
\end{cases}
\end{align*}
where
\[
a=\frac{\lambda}{n}\min\{B^1(n-B^1-G^1),B^2(n-B^2-G^2)\}.
\]
By direct calculation, it is easy to check that the Markov process
with the above transition rates function is a coupling of
$(B_t^1,G_t^1)$ and $(B_t^2,G_t^2)$ while all the state transitions
hold the property that $(B^1,G^1)\succeq (B^2,G^2)$.

\qed

\section{Mean field limit}\label{section three}
In this section we introduce an ODE, the solution to which is the
mean field limit of $(\frac{B_t}{n},\frac{G_t}{n})$ as $n$ grows to
infinity. For later use, for any $x=(b,g)\in \mathbb{R}^2$, we use
$\|x\|_1$ to denote $|b|+|g|$, which is the $l_1$ norm of $x=(b,g)$.

We consider the following two-dimensional ODE with initial condition
$(b_0,g_0)\in [0,1]\times[0,1]$.
\begin{equation*}
\begin{cases}
\frac{d}{dt}b_t=-b_t+\lambda b_tg_t,\\
\frac{d}{dt}g_t=-g_t-\lambda b_tg_t+\lambda b_t(1-b_t-g_t),
\end{cases}
\end{equation*}
where $\lambda>0$. For simplicity, we write the above ODE as
\begin{equation}\label{equ ODE 3.1}
\frac{d}{dt}(b_t,g_t)=F(b_t,g_t),
\end{equation}
where $F=(F_1,F_2)$ such that $F_1(b,g)=-b+\lambda bg$ and
$F_2(b,g)=-g-\lambda bg+\lambda b(1-b-g)$.

Later we will show that $\frac{1}{n}(B_t,G_t)$ converges to the
solution $(b_t,g_t)$ to ODE \eqref{equ ODE 3.1} in probability. By
direct calculation, it is easy to check that
\[
F(b,g)=(0,0)
\]
has the unique solution $(0,0)$ when $\lambda<4$ while has two
solutions in $(0,1)\times(0,1)$ when $\lambda>4$, which intuitively
explains why the critical value of our process is $4$.

We define $\Lambda=\{(b,g):~b\geq 0,g\geq 0,b+g\leq 1\}\subseteq
\mathbb{R}^2$, then we have the following lemma.
\begin{lemma}\label{lemma 3.1}
The solution $(b_t,g_t)$ to ODE \eqref{equ ODE 3.1} with initial
condition $(b_0,g_0)\in \Lambda$ exists for $t\in [0,+\infty)$ and
is unique. Furthermore, for any $t\geq 0$, $(b_t,g_t)\in \Lambda$.
\end{lemma}
\proof

It is easy to check that $F$ satisfies the local Lipschitz condition
under the norm $\|\cdot\|_1$, according to which the uniqueness of
the solution holds. On the boundary of $\Lambda$, it is easy to
check that any vector of the vector-field of ODE \eqref{equ ODE 3.1}
points to the inner of $\Lambda$, hence the solution is absorbed in
the area $\Lambda$. On the area $\Lambda$, it is easy to check that
$F$ satisfies the global Lipschitz condition, hence the solution
exists for $t\in [0,+\infty)$.

\qed

The next lemma shows that $(\frac{B_t}{n},\frac{G_t}{n})$ converges
to the solution $(b_t,g_t)$ to the ODE \eqref{equ ODE 3.1} in
probability as $n$ grows to infinity.
\begin{lemma}\label{Lemma 3.2}
Let $\{(b_t,g_t)\}_{t\geq 0}$ be the solution to the ODE \eqref{equ
ODE 3.1} with initial condition $(b_0,g_0)\in \Lambda$, then for any
$T>0$ and $\epsilon>0$, there exist constants $C_1=C_1(T,\epsilon)$
and $N_1=N_1(T,\epsilon)$ such that
\[
P_{\lambda,n}^{nb_0,ng_0}\Big(\sup_{0\leq t\leq
T}\big\|(\frac{B_t}{n},\frac{G_t}{n})-(b_t,g_t)\big\|_1~>\epsilon\Big)\leq
\frac{C_1}{n}
\]
for any $n\geq N_1$.
\end{lemma}
Note that $C_1$ and $N_1$ do not depend on the choice of
$(b_0,g_0)$.

The proof of Lemma \ref{Lemma 3.2} follows the analysis introduced
by Ethier and Kurtz to construct the theory of density-dependent
population model (See Chapter 11 of \cite{Ethier1986}). Readers
familiar with this theory can skip the following proof.

\proof[Proof of Lemma \ref{Lemma 3.2}]

By the definition of $(b_t,g_t)$,
\begin{equation}\label{equ 3.2}
\begin{cases}
b_t=b_0+\int_0^tF_1(b_s,g_s)~ds,\\
g_t=g_0+\int_0^tF_2(b_s,g_s)~ds.
\end{cases}
\end{equation}
Let $\{N_j(t):t\geq 0\}_{j=1,2,3,4}$ be four independent copies of
the Poisson process with rate one, then according to the transition
rates function of $\{B_t,G_t\}_{t\geq 0}$ and Theorem 6.4.1 of
\cite{Ethier1986}, we can write $(B_t,G_t)$ as
\begin{equation}\label{equ 3.3}
\begin{cases}
B_t=\lfloor nb_0\rfloor-N_1(\int_0^tB_s~d_s)+N_3(\int_0^t\frac{\lambda}{n}B_sG_s~ds),\\
G_t=\lfloor
ng_0\rfloor-N_2(\int_0^tG_s~ds)-N_3(\int_0^t\frac{\lambda}{n}B_sG_s~ds)+N_4(\int_0^t\frac{\lambda}{n}B_s(n-B_s-G_s)~ds).
\end{cases}
\end{equation}
For $j=1,2,3,4$, we define $\widetilde{N}_j(t)=N_j(t)-t$, then
$\{\widetilde{N}_j(t):~t\geq 0\}$ is a martingale with
$E\widetilde{N}_j(t)=0$ and $E[\widetilde{N}_j^2(t)]=t$. According
to Equations \eqref{equ 3.2}, \eqref{equ 3.3} and the definition of
$F=(F_1,F_2)$,
\begin{equation}\label{equ 3.4}
\begin{cases}
\frac{B_t}{n}-b_t=&\frac{\lfloor
nb_0\rfloor-nb_0}{n}+\int_0^tF_1(\frac{B_s}{n},\frac{G_s}{n})-F_1(b_s,g_s)~ds-\frac{\widetilde{N}_1(\int_0^tB_s~d_s)}{n}
+\frac{\widetilde{N}_3(\int_0^t\frac{\lambda}{n}B_sG_s~ds)}{n},\\
\frac{G_t}{n}-g_t=&\frac{\lfloor ng_0\rfloor-ng_0}{n}+\int_0^tF_2(\frac{B_s}{n},\frac{G_s}{n})-F_2(b_s,g_s)~ds\\
&-\frac{\widetilde{N}_2(\int_0^tG_s~ds)}{n}
-\frac{\widetilde{N}_3(\int_0^t\frac{\lambda}{n}B_sG_s~ds)}{n}+\frac{\widetilde{N}_4(\int_0^t\frac{\lambda}{n}B_s(n-B_s-G_s)~ds)}{n}.
\end{cases}
\end{equation}
Since $B_t\leq n$, $|\widetilde{N}_1(\int_0^tB_s~d_s)|\leq
\max\limits_{0\leq s\leq nt}|\widetilde{N}_1(s)|$. For similar
reasons,
\begin{align*}
|\widetilde{N}_2(\int_0^tB_s~d_s)|&\leq \max\limits_{0\leq s\leq
nt}|\widetilde{N}_1(s)|,\\
|\widetilde{N}_3(\int_0^t\frac{\lambda}{n}B_sG_s~ds)|&\leq
\max\limits_{0\leq s\leq \lambda nt}|\widetilde{N}_3(s)|,\\
|\widetilde{N}_4(\int_0^t\frac{\lambda}{n}B_sG_s~ds)|&\leq
\max\limits_{0\leq s\leq \lambda nt}|\widetilde{N}_4(s)|.
\end{align*}
It is easy to check that $F=(F_1,F_2)$ satisfies the global
Lipschitz condition on $\Lambda$ under the norm $\|\cdot\|_1$, hence
there exists $K>0$ such that
\[
\|F(b_1,g_1)-F(b_2,g_2)\|_1\leq K\|(b_1,g_1)-(b_2,g_2)\|_1
\]
for any $(b_1,g_1),(b_2,g_2)\in \Lambda$. It is obviously that
$\frac{1}{n}(B_t,G_t)\in \Lambda$ for ant $t\geq 0$. As we have
shown in Lemma \ref{lemma 3.1}, $(b_t,g_t)\in \Lambda$ for any
$t\geq 0$. As a result, by Equation \eqref{equ 3.4},
\begin{equation}\label{equ 3.5}
\big\|(\frac{B_t}{n},\frac{G_t}{n})-(b_t,g_t)\big\|_1 \leq
\int_0^tK\big\|(\frac{B_s}{n},\frac{G_s}{n})-(b_s,g_s)\big\|_1~ds+\frac{M(\lambda,n,t)}{n}
\end{equation}
for any $t\geq 0$, where
\[
M(\lambda,n,t)=2+\max_{0\leq s\leq
nt}|\widetilde{N}_1(s)|+\max_{0\leq s\leq
nt}|\widetilde{N}_2(s)|+2\max_{0\leq s\leq \lambda
nt}|\widetilde{N}_3(s)|+\max_{0\leq s\leq \lambda
nt}|\widetilde{N}_4(s)|.
\]
Since $M(\lambda,n,t)$ increases with $t$,
\[
\big\|(\frac{B_t}{n},\frac{G_t}{n})-(b_t,g_t)\big\|_1 \leq
\int_0^tK\big\|(\frac{B_s}{n},\frac{G_s}{n})-(b_s,g_s)\big\|_1~ds+\frac{M(\lambda,n,T)}{n}
\]
for $0\leq t\leq T$. Then by Grownwall's inequality,
\begin{equation}\label{equ 3.6}
\big\|(\frac{B_t}{n},\frac{G_t}{n})-(b_t,g_t)\big\|_1\leq
\frac{M(\lambda,n,T)}{n}e^{Kt}
\end{equation}
for $0\leq t\leq T$. By Equation \eqref{equ 3.6},
\begin{equation}\label{equ 3.6 two}
\sup_{0\leq t\leq
T}\big\|(\frac{B_t}{n},\frac{G_t}{n})-(b_t,g_t)\big\|_1~\leq
\frac{M(\lambda,n,T)}{n}e^{KT}.
\end{equation}
Since $\{\widetilde{N}_1(t):~t\geq 0\}$ is a martingale, according
to Doob's inequality,
\[
P\Big(\max_{0\leq s\leq nT}|\widetilde{N}_1(s)|\geq n
\epsilon_0\Big)\leq
\frac{1}{n^2\epsilon_0^2}E[\widetilde{N}^2_1(nT)]=\frac{T}{n\epsilon_0^2}
\]
for any $\epsilon_0>0$. For given $\epsilon>0$, we choose
$\epsilon_0=\frac{\epsilon}{8e^{KT}}$, then
\[
P\Big(\max_{0\leq s\leq nT}|\widetilde{N}_1(s)|\geq
\frac{n\epsilon}{8e^{KT}}\Big)\leq \frac{64e^{2KT}T}{n\epsilon^2}.
\]
According similar analysis,
\begin{align*}
P\Big(\max_{0\leq s\leq nT}|\widetilde{N}_2(s)|\geq
\frac{n\epsilon}{8e^{KT}}\Big)&\leq
\frac{64e^{2KT}T}{n\epsilon^2},\\
P\Big(2\max_{0\leq s\leq \lambda nT}|\widetilde{N}_3(s)|\geq
\frac{n\epsilon}{8e^{KT}}\Big)&\leq \frac{256e^{2KT}\lambda
T}{n\epsilon^2},\\
P\Big(\max_{0\leq s\leq \lambda nT}|\widetilde{N}_4(s)|\geq
\frac{n\epsilon}{8e^{KT}}\Big)&\leq \frac{64e^{2KT}\lambda
T}{n\epsilon^2}.
\end{align*}
Therefore,
\[
P\Big(M(\lambda,n,T)\geq2+\frac{n\epsilon}{2e^{KT}}\Big)\leq
\frac{128e^{2KT}T+320e^{2KT}\lambda T}{\epsilon^2}.
\]
As a result,
\[
P\big(\frac{M(\lambda,n,T)}{n}\geq \frac{\epsilon}{e^{KT}}\big)\leq
\frac{C_1(T,\epsilon)}{n}
\]
for $n\geq N_1(T,\epsilon)=\frac{4e^{KT}}{\epsilon}$, where
\[
C_1(T,\epsilon)=\frac{128e^{2KT}T+320e^{2KT}\lambda T}{\epsilon^2}.
\]
Then by Equation \eqref{equ 3.6 two},
\[
P_{\lambda,n}^{n,0}\Big(\sup_{0\leq t\leq
T}\big\|(\frac{B_t}{n},\frac{G_t}{n})-(b_t,g_t)\big\|_1~>\epsilon\Big)
\leq P\big(\frac{M(\lambda,n,T)}{n}\geq
\frac{\epsilon}{e^{KT}}\big)\leq \frac{C_1(T,\epsilon)}{n}
\]
for $n\geq N_1(T,\epsilon)$ and the proof is complete.

\qed

The next lemma is crucial for the proof of Equation \eqref{equ 2.2
main}.
\begin{lemma}\label{lemma 3.3}
Let $\{(b_t,g_t)\}_{t\geq 0}$ be the solution to ODE \eqref{equ ODE
3.1} with initial condition $(b_0,g_0)=(1,0)$, then when
$\lambda<4$,
\[
\lim_{t\rightarrow+\infty}\|(b_t,g_t)\|_1=0.
\]
\end{lemma}

\proof

It is easy to check that $F_2(b,g)\geq 0$ when and only when $g\leq
\frac{\lambda b(1-b)}{2\lambda b+1}$. Let $g^*=\max\{\frac{\lambda
b(1-b)}{2\lambda b+1}:~0\leq b\leq 1\}$, then it is easy to check
that $g^*<\frac{1}{\lambda}$ when $\lambda<4$. Let
$\widetilde{g}=\frac{\frac{1}{\lambda}+g^*}{2}$, then
$g_0=0<\widetilde{g}$. For any $0\leq b\leq 1$,
$F_2(b,\widetilde{g})<0$ since $\widetilde{g}>g^*$. As a result,
$g_t$ can never exceed $\widetilde{g}$ since
$\frac{d}{dt}g_t\Big|_{g_t=\widetilde{g}}=F_2(b_t,\widetilde{g})<0$.
Therefore, $g_t\leq \widetilde{g}$ for any $t\geq 0$. Then,
\[
\frac{d}{dt}b_t\leq -b_t+\lambda \widetilde{g}b_t=(\lambda
\widetilde{g}-1)b_t
\]
and hence
\begin{equation}\label{equ 3.7}
b_t\leq e^{(\lambda \widetilde{g}-1)t}
\end{equation}
for any $t\geq 0$. By Equation \eqref{equ 3.7},
\[
\frac{d}{dt}g_t=F_2(b_t,g_t)\leq -g_t+\lambda e^{(\lambda
\widetilde{g}-1)t}
\]
for any $t\geq 0$. As a result,
\[
\frac{d}{dt}(e^tg_t)\leq \lambda e^{\lambda\widetilde{g}t}
\]
and hence
\begin{equation}\label{equ 3.8}
g_t\leq \frac{e^{(\lambda\widetilde{g}-1)t}-e^{-t}}{\widetilde{g}}.
\end{equation}
Since $\widetilde{g}<\frac{1}{\lambda}$, $\lambda
\widetilde{g}-1<0$. As a result, Lemma \ref{lemma 3.3} follows from
Equations \eqref{equ 3.7} and \eqref{equ 3.8} directly.

\qed

\section{Proof of Equation \eqref{equ 2.1 main}}\label{section four}
In this section we give the proof of Equation \eqref{equ 2.1 main}.
The intuitive idea of the proof is as follows. When $\lambda>4$, it
is easy to check that there exists $(b_0,g_0)$ in the inner of
$\Lambda$ such that $F_1(b_0,g_0)>0$ and $F_2(b_0,g_0)>0$. Then, by
analyzing the vector-field of ODE \eqref{equ ODE 3.1}, it is easy to
check that the solution $\{(b_t,g_t)\}_{t\geq 0}$ with initial
condition $(b_0,g_0)$ is absorbed in the area
\[
\Lambda_1=\{(b,g):~b\geq b_0,b+g\geq b_0+g_0\}.
\]
As shown in Lemma \ref{Lemma 3.2}, conditioned on
$(B_0,G_0)=(nb_0,ng_0)$, $\frac{1}{n}(B_t,G_t)$ is approximate to
$(b_t,g_t)$. Hence $(B_t,G_t)$ should stay in the area $n\Lambda_1$
for a long time. Since $(nb_0,ng_0)\preceq (n,0)$ and
$\{(B_t,G_t)\}_{t\geq 0}$ is monotone under the partial order
$\preceq$ by Lemma \ref{lemma 2.2}, $(B_t,G_t)$ with
$(B_0,G_0)=(n,0)$ should stay in $n\Lambda$ for a longer time.

Our proof is the effort to make the above intuitive idea rigorous
and show that the precise meaning of `long time' is $\exp\{O(n)\}$
units of time. First we show the existence of $(b_0,g_0)$. We define
\[
\widetilde{\Lambda}=\{(b,g):~b>0,g>0,b+g<1\}
\]
as the inner of $\Lambda$.
\begin{lemma}\label{lemma 4.1}
When $\lambda>4$, there exists $(b_0,g_0)\in \widetilde{\Lambda}$
such that
\[
F_1(b_0,g_0)>0\text{~and~}F_2(b_0,g_0)>0,
\]
where $F=(F_1,F_2)$ is as defined in ODE \eqref{equ ODE 3.1}.
\end{lemma}

\proof

By direct calculation, when $\lambda>4$,
$(\frac{\lambda-2}{2\lambda},\frac{1}{\lambda})\in
\widetilde{\Lambda}$ and
\[
F_2(\frac{\lambda-2}{2\lambda},\frac{1}{\lambda})>0.
\]
Then, we can choose sufficiently small positive $\beta$ such that
$F_2(\frac{\lambda-2}{2\lambda},\frac{1}{\lambda}+\beta)>0$ and
$(\frac{\lambda-2}{2\lambda},\frac{1}{\lambda}+\beta)\in
\widetilde{\Lambda}$. For any $(b,g)\in \Lambda$ with $g>
\frac{1}{\lambda}$, $F_1(b,g)>0$. As a result, Lemma \ref{lemma 4.1}
holds with
\[
(b_0,g_0)=(\frac{\lambda-2}{2\lambda},\frac{1}{\lambda}+\beta).
\]
\qed

For later use, we choose $\alpha>0$ sufficiently small such that
$(1-\alpha)(b_0,g_0),(1+\alpha)(b_0,g_0)\in \widetilde{\Lambda}$,
$\overline{g},\underline{g}>0$ and
\begin{align}\label{equ 4.1}
&\lambda(1-\alpha)b_0\underline{g}>(1+\alpha)b_0, \\
&\lambda
(1-\alpha)b_0[1-(1+\alpha)(b_0+g_0)]>\overline{g}+\lambda(1+\alpha)b_0\overline{g},
\notag
\end{align}
where $\underline{g}=(1-\alpha)(b_0+g_0)-(1+\alpha)b_0$ and
$\overline{g}=(1+\alpha)(b_0+g_0)-(1-\alpha)b_0$. Note that the
existence of $\alpha$ depends on the fact that
\begin{align*}
&\lambda b_0g_0>b_0, \\
&\lambda b_0[1-(b_0+g_0)]>g_0+\lambda b_0g_0.
\end{align*}
according to Lemma \ref{lemma 4.1}.

We define
\[
\Lambda_2=\{(b,g):~(1-\alpha)b_0\leq b\leq (1+\alpha)b_0,
(1-\alpha)(b_0+g_0)\leq b+g\leq (1+\alpha)(b_0+g_0)\},
\]
then it is easy to check that
\begin{equation}\label{equ 4.1 two}
\underline{g}=\min\{g:~(b,g)\in
\Lambda_2\}\text{~and~}\overline{g}=\max\{g:~(b,g)\in \Lambda_2\}.
\end{equation}
We define
\[
\gamma=\inf\{t\geq 0: ~\frac{1}{n}(B_t,G_t)\not \in \Lambda_2\}
\]
as the first moment that $(B,G)$ exits $n\Lambda_2$. Furthermore, we
define
\[
D=\inf\{\|(b,g)-(b_0,g_0)\|_1:~(b,g)\not\in \Lambda_2\}>0.
\]
The next lemma about the time $(B,G)$ waits for to exit $n\Lambda_2$
conditioned on $(B_0,G_0)=(nb_0,ng_0)$ is utilized later.
\begin{lemma}\label{lemma 4.2}
For given $\lambda>4$ and $\alpha$ defined as in Equation
\eqref{lemma 4.1}, there exists $C_4=C_4(\lambda,\alpha)>0$ and
$N_6=N_6(\lambda,\alpha)>0$ such that
\[
P^{nb_0,ng_0}_{\lambda,n}\Big(\gamma\geq \frac{
D}{4(1+\lambda)}\Big)\geq 1-e^{-C_4n}
\]
for $n\geq N_6$.
\end{lemma}
\proof

Conditioned on $(B_0,G_0)=(nb_0,ng_0)$, the $l_1$ norm $\|(B,G)\|_1$
must change by at least $nD$ for $\frac{1}{n}(B,G)$ to exit
$\Lambda_2$. At each moment that $(B,G)$ jumps, $\|(B,G)\|_1$
changes by at most $2$. Hence $(B,G)$ must jump at least
$\frac{nD}{2}$ times to exit $n\Lambda_2$. It is easy to check that
$(B,G)$ changes state with rate at most
\[
n+\frac{\lambda}{n}nn=(1+\lambda)n.
\]
As a result,
\begin{equation*}
P^{nb_0,ng_0}_{\lambda,n}\Big(\gamma\leq
\frac{D}{4(1+\lambda)}\Big)\leq
P\Big(Y_n\big(\frac{D}{4(1+\lambda)}\big)\geq \frac{nD}{2}\Big),
\end{equation*}
where $\{Y_n(t)\}_{t\geq 0}$ is the Poisson process with rate
$(1+\lambda)n$ for each $n\geq 1$. $Y_n(t)$ has the same probability
distribution as that of $Y_1(nt)$, hence
\begin{equation}\label{equ 4.2}
P^{nb_0,ng_0}_{\lambda,n}\Big(\gamma\leq
\frac{D}{4(1+\lambda)}\Big)\leq
P\Big(\frac{Y_1\big(n\frac{D}{4(1+\lambda)}\big)}{n}\geq
\frac{D}{2}\Big).
\end{equation}
According to classic limit theorems of the Poisson process, for any
$t>0$, $Y_1(nt)/n$ converges to $(1+\lambda)t$ in probability as
$n\rightarrow+\infty$ and there exists $I(t)>0$ such that
\[
P\Big(\frac{Y_1(nt)}{n}\geq 2(1+\lambda)t\Big)\leq e^{-nI(t)}
\]
for sufficiently large $n$. As a result, there exists
$N_6=N_6(\lambda,\alpha)$ such that
\begin{equation}\label{equ 4.3}
P\Big(\frac{Y_1\big(n\frac{D}{4(1+\lambda)}\big)}{n}\geq
\frac{D}{2}\Big)\leq e^{-nI(\frac{D}{4(1+\lambda)})}
\end{equation}
for $n\geq N_6$. Lemma \ref{lemma 4.2} follows from Equations
\eqref{equ 4.2} and \eqref{equ 4.3} directly with
$C_4=I(\frac{D}{4(1+\lambda)})$.

\qed

We introduce a birth-and-death process as an auxiliary model for the
proof of Equation \eqref{equ 2.1 main}. Let
$\{(\widehat{B}_t,\widehat{S}_t)\}_{t\geq 0}$ be the birth-and-death
process with transition rates function given by
\begin{equation*}
(\widehat{B},\widehat{S})\rightarrow
\begin{cases}
(\widehat{B},\widehat{S})-(1,1) &\text{~at rate~} n(1+\alpha)b_0,\\
(\widehat{B},\widehat{S})-(0,1) &\text{~at rate~} n\overline{g},\\
(\widehat{B},\widehat{S})+(1,0)  &\text{~at rate~} \lambda
n(1-\alpha)b_0\underline{g},\\
(\widehat{B},\widehat{S})+(0,1) &\text{~at rate~} \lambda
n(1-\alpha)b_0[1-(1+\alpha)(b_0+g_0)].
\end{cases}
\end{equation*}
Later we will show that $B_t\geq \widehat{B}_t$ and $B_t+G_t\geq
\widehat{S}_t$ in the sense of coupling for $t\in [0,\gamma]$, for
which we introduce $(\widehat{B},\widehat{S})$. The next lemma shows
that for any $t>0$, $\widehat{S}_t\geq \widehat{S}_0$ and
$\widehat{B}_t\geq \widehat{B}_0$ with high probability.
\begin{lemma}\label{lemma 4.3}
There exists $C_5=C_5(\lambda,\alpha)>0$ such that for any $t>0$ and
$n\geq 1$,
\[
P\Big(\widehat{B}_t\geq nb_0,\widehat{S}_t\geq
ng_0+nb_0\Big|(\widehat{B}_0,\widehat{S}_0)=\big(nb_0,ng_0+nb_0\big)\Big)\geq
1-2e^{-C_5nt}.
\]
\end{lemma}

\proof Throughout this proof we assume that
$(\widehat{B}_0,\widehat{S}_0)=(nb_0,nb_0+ng_0)$. Let $q_1=\lambda
(1-\alpha)b_0[1-(1+\alpha)(b_0+g_0)]$ and
$q_2=(1+\alpha)b_0+\overline{g}$, then $\{\widehat{S}_t\}_{t\geq 0}$
is a birth-and-death process with transition rates function given by
\[
\widehat{S}\rightarrow
\begin{cases}
\widehat{S}+1 &\text{~at rate~} nq_1,\\
\widehat{S}-1 &\text{~at rate~} nq_2.
\end{cases}
\]
By Equation \eqref{equ 4.1}, it is easy to check that $q_1>q_2$,
then we can choose $\frac{q_2}{q_1}<\rho<1$. Since $\rho<1$, by
Chebyshev's inequality,
\begin{equation}\label{equ 4.5}
P(\widehat{S}_t\leq nb_0+ng_0)=P(\rho^{\widehat{S}_t}\geq
\rho^{nb_0+ng_0})\leq \rho^{-nb_0-ng_0}E\rho^{\widehat{S}_t}.
\end{equation}
According to the transition rates function of $\widehat{S}$,
\begin{align*}
\frac{d}{dt}E\rho^{\widehat{S}_t}&=nq_1\big(E\rho^{\widehat{S}_t+1}-E\rho^{\widehat{S}_t}\big)
+nq_2\big(E\rho^{\widehat{S}_t-1}-E\rho^{\widehat{S}_t}\big)\\
&=n\big(q_1\rho+\frac{q_2}{\rho}-q_1-q_2\big)E\rho^{\widehat{S}_t}.
\end{align*}
Then,
\begin{equation}\label{equ 4.6}
E\rho^{\widehat{S}_t}=\rho^{nb_0+ng_0}e^{-C_6nt}
\end{equation}
since $\widehat{S}_0=nb_0+ng_0$, where
$C_6=q_1+q_2-q_1\rho-\frac{q_2}{\rho}$. Note that $C_6>0$ since
$\frac{q_2}{q_1}<\rho<1$. By Equations \eqref{equ 4.5} and
\eqref{equ 4.6},
\begin{equation}\label{equ 4.7}
P(\widehat{S}_t\leq nb_0+ng_0)\leq e^{-C_6nt}.
\end{equation}
According to the same analysis as that gives Equation \eqref{equ
4.7}, there exists $C_7>0$ such that
\begin{equation}\label{equ 4.8}
P(\widehat{B}_t\leq nb_0)\leq e^{-C_7nt}.
\end{equation}
Let $C_5=\min\{C_6,C_7\}$, then Lemma \ref{lemma 4.3} follows from
Equations \eqref{equ 4.7} and \eqref{equ 4.8} directly.

\qed

The next lemma shows that $B_t$ and $B_t+G_t$ are bounded from below
by $\widehat{B}_t$ and $\widehat{S}_t$ respectively for $t\geq
[0,\gamma]$.
\begin{lemma}\label{lemma 4.4}
Conditioned on
$(B_0,B_0+G_0)=(\widehat{B}_0,\widehat{S}_0)=(nb_0,nb_0+ng_0)$,
\[
B_t\geq \widehat{B}_t \text{~and~}B_t+G_t\geq \widehat{S}_t
\]
in the sense of coupling for $t\in [0,\gamma]$.
\end{lemma}

\proof

For $t\in [0,\gamma]$, $\frac{1}{n}(B_t,G_t)\in \Lambda_2$. The
transition rates function of $\{(B_t,B_t+G_t)\}_{t\geq 0}$ is given
by
\[
(B,B+G)\rightarrow
\begin{cases}
(B,B+G)-(1,1)&\text{~at rate~} F_1(B,G)=B,\\
(B,B+G)-(0,1)&\text{~at rate~} F_2(B,G)=G,\\
(B,B+G)+(1,0)&\text{~at rate~} F_3(B,G)=\frac{\lambda}{n}BG,\\
(B,B+G)+(0,1)&\text{~at rate~} F_4(B,G)=\frac{\lambda}{n}B(n-B-G).
\end{cases}
\]
According to the definition of $\Lambda_2$, for any $(B,G)$ that
$\frac{1}{n}(B,G)\in \Lambda_2$,
\begin{align*}
&F_1(B,G)\leq n(1+\alpha)b_0, \text{~~}F_2(B,G)\leq n\overline{g},\\
&F_3(B,G)\geq\lambda
n(1-\alpha)b_0\underline{g},\text{~~}F_4(B,G)\geq\lambda
n(1-\alpha)b_0[1-(1+\alpha)(b_0+g_0)].
\end{align*}
As a result, we can couple $(B_t,B_t+G_t)$ and
$(\widehat{B}_t,\widehat{S}_t)$ as follows for $t\in [0,\gamma]$.
\begin{align*}
&(B,B+G,\widehat{B},\widehat{S})\\
&\rightarrow
\begin{cases}
(B,B+G,\widehat{B},\widehat{S})-(1,1,1,1)&\text{~at rate~} F_1(B,G),\\
(B,B+G,\widehat{B},\widehat{S})-(0,0,1,1)&\text{~at rate~}
n(1+\alpha)b_0-F_1(B,G),\\
(B,B+G,\widehat{B},\widehat{S})-(0,1,0,1)&\text{~at rate~} F_2(B,G),\\
(B,B+G,\widehat{B},\widehat{S})-(0,0,0,1)&\text{~at rate~}
n\overline{g}-F_2(B,G),\\
(B,B+G,\widehat{B},\widehat{S})+(1,0,1,0)&\text{~at rate~} \lambda
n(1-\alpha)b_0\underline{g},\\
(B,B+G,\widehat{B},\widehat{S})+(1,0,0,0)&\text{~at rate~}
F_3(B,G)-\lambda n(1-\alpha)b_0\underline{g},\\
(B,B+G,\widehat{B},\widehat{S})+(0,1,0,1)&\text{~at rate~} \lambda
n(1-\alpha)b_0[1-(1+\alpha)(b_0+g_0)],\\
(B,B+G,\widehat{B},\widehat{S})+(0,1,0,0)&\text{~at rate~}
F_4(B,G)-\lambda n(1-\alpha)b_0[1-(1+\alpha)(b_0+g_0)].
\end{cases}
\end{align*}
The above coupling does not change the property that $B\geq
\widehat{B}$ and $B+G\geq \widehat{S}$, hence the proof is complete.

\qed

At last we give the proof of Equation \eqref{equ 2.1 main}.

\proof[Proof of Equation \eqref{equ 2.1 main}]

We define
\[
\Lambda_1=\{(b,g)\in \Lambda:~b\geq b_0, b+g\geq b_0+g_0\}
\]
as we have done at the beginning of this section.

Let $T_4=\frac{D}{4(1+\lambda)}>0$, where $D$ is defined as in Lemma
\ref{lemma 4.2}, then the first step of this proof is to show that
\begin{equation}\label{equ 4.9}
P_{\lambda,n}^{nb_0,ng_0}\big(\frac{(B_{T_4},G_{T_4})}{n}\in
\Lambda_1\big)\geq 1-3e^{-C_8n}
\end{equation}
for some $C_8=C_8(\lambda)>0$ and sufficiently large $n$. The proof
of Equation \eqref{equ 4.9} is as follows. By Lemma \ref{lemma 4.4},
$B_t\geq \widehat{B}_t$ and $B_t+G_t\geq \widehat{S}_t$ for $t\in
[0,\gamma]$, where
$(B_0,G_0+B_0)=(\widehat{B}_0,\widehat{S}_0)=(nb_0,nb_0+ng_0)$.
Therefore,
\begin{align}\label{equ 4.10}
P_{\lambda,n}^{nb_0,ng_0}\big(\frac{(B_{T_4},G_{T_4})}{n}\in
\Lambda_1\big) &\geq
P_{\lambda,n}^{nb_0,ng_0}\big(\frac{(B_{T_4},G_{T_4})}{n}\in
\Lambda_1,~\gamma\geq T_4\big) \\
&=P_{\lambda,n}^{nb_0,ng_0}\big(B_{T_4}\geq nb_0,B_{T_4}+G_{T_4}\geq nb_0+ng_0,~\gamma\geq T_4\big)\notag\\
&\geq P_{\lambda,n}^{nb_0,ng_0}\big(\widehat{B}_{T_4}\geq
nb_0,\widehat{S}_{T_4}\geq nb_0+ng_0,~\gamma\geq T_4\big)\notag\\
&\geq P(\widehat{B}_{T_4}\geq nb_0,\widehat{S}_{T_4}\geq
nb_0+ng_0)-P_{\lambda,n}^{nb_0,ng_0}\big(\gamma<T_4\big)\notag\\
&\geq 1-2e^{-C_5T_4n}-e^{-C_4n}\notag
\end{align}
for $n\geq N_6$, where $N_6,C_4$ are defined as in Lemma \ref{lemma
4.2} while $C_5$ is defined as in Lemma \ref{lemma 4.3}. Note that
the last inequality in Equation \eqref{equ 4.10} follows from Lemmas
\ref{lemma 4.2} and \ref{lemma 4.3}. As a result, Equation
\eqref{equ 4.9} holds with $C_8=\min\{C_5T_4,C_4\}$ and $n\geq N_6$.

Note that if $(B_1,G_1)\succeq (B_2,G_2)$ and
$\frac{1}{n}(B_2,G_2)\in \Lambda_1$, then $\frac{1}{n}(B_1,G_1)\in
\Lambda_1$. Therefore, by Lemma \ref{lemma 2.2},
\[
P_{\lambda,n}^{nb_0,ng_0}\big(\frac{(B_{T_4},G_{T_4})}{n}\in
\Lambda_1\big)\leq
P_{\lambda,n}^{B,G}\big(\frac{(B_{T_4},G_{T_4})}{n}\in
\Lambda_1\big)
\]
for any $(B,G)\succeq (nb_0,ng_0)$. Then by Equation \eqref{equ
4.9},
\begin{equation}\label{equ 4.11}
P_{\lambda,n}^{B,G}\big((B_{T_4},G_{T_4})\succeq (nb_0,ng_0)\big)=
P_{\lambda,n}^{B,G}\big(\frac{(B_{T_4},G_{T_4})}{n}\in
\Lambda_1\big)\geq 1-3e^{-C_8n}
\end{equation}
for any $(B,G)\succeq (nb_0,ng_0)$ and $n\geq N_6$. By Equation
\eqref{equ 4.11}, utilizing the Markov property for
$e^{\frac{C_8}{2}n}$ times,
\begin{equation*}
P_{\lambda,n}^{nb_0,ng_0}\Big(\big(B_{T_4e^{\frac{C_8}{2}n}},G_{T_4e^{\frac{C_8}{2}n}}\big)\succeq
(nb_0,ng_0)\Big)\geq (1-3e^{-C_8n})^{e^{\frac{C_8}{2}n}}\geq
1-3e^{-\frac{C_8}{2}n}
\end{equation*}
for $n\geq N_6$ and hence
\begin{equation}\label{equ 4.12}
P_{\lambda,n}^{nb_0,ng_0}(\tau>T_4e^{\frac{C_8}{2}n})\geq
1-3e^{-\frac{C_8}{2}n}
\end{equation}
for $n\geq N_6$. By Equation \eqref{equ 4.12},
\begin{equation}\label{equ 4.13}
P_{\lambda,n}^{nb_0,ng_0}(\tau>e^{\frac{C_8}{4}n})\geq
1-3e^{-\frac{C_8}{2}n}
\end{equation}
for sufficiently large $n$. Then, by Lemma \ref{lemma 2.2} and the
fact that $(n,0)\succeq (nb_0,ng_0)$,
\[
P_{\lambda,n}^{n,0}(\tau>e^{\frac{C_8}{4}n})\geq
P_{\lambda,n}^{nb_0,ng_0}(\tau>e^{\frac{C_8}{4}n})\geq
1-3e^{-\frac{C_8}{2}n}
\]
for sufficiently large $n$ and hence Equation \eqref{equ 2.1 main}
holds with $C(\lambda)=\frac{C_8}{4}$.

\qed

\section{Proof of Equation \eqref{equ 2.2 main}}\label{section five}
In this section we give the proof of Equation \eqref{equ 2.2 main}.
It is obviously that we only need to deal with small $\theta$, so we
assume that $\theta<1$. Throughout this section we use
$\{(b_t,g_t)\}_{t\geq 0}$ to denote the solution to ODE \eqref{equ
ODE 3.1} with initial condition $(b_0,g_0)=(1,0)$. Sometimes we
write $b_t,g_t$ as $b(t),g(t)$ when the subscript is complex. First
we show that $\frac{1}{n}\|(B_t,G_t)\|_1$ stays small for
$O(\sqrt{n})$ units of time with high probability when $\lambda<4$.
\begin{lemma}\label{lemma 5.1}
For given $\lambda<4$ and $\theta\in(0,1)$, there exist
$T_1=T_1(\lambda,\theta)$, $N_2=N_2(\lambda,\theta)$,
$C_2=C_2(\lambda,\theta)$ and $C_3=C_3(\lambda,\theta)$ such that
\[
P_{\lambda,n}^{n,0}\Big(\sup_{0\leq t\leq
C_2\sqrt{n}}\frac{\|(B_{T_1+t},G_{T_1+t})\|_1}{n}~\geq
\frac{\theta}{(3+\theta)\lambda}\Big)\leq \frac{C_3}{\sqrt{n}}
\]
for each $n\geq N_2$.
\end{lemma}

\proof

By Lemma \ref{lemma 3.3}, for $\lambda<4$, we can choose
$0<T_3(\lambda,\theta)<T_1(\lambda,\theta)$ such that
\[
b(T_3)+g(T_3)\leq \frac{\theta}{2(3+\theta)\lambda},
\text{~}b(T_1)+g(T_1)\leq \frac{b(T_3)}{2}
\]
and
\[
b_t+g_t\leq \frac{\theta}{2(3+\theta)\lambda}
\]
for any $t\in [T_3, T_1]$. Let
$\{(\widehat{b}_t,\widehat{g}_t)\}_{t\geq 0}$ be the solution to ODE
\eqref{equ ODE 3.1} with
$(\widehat{b}_0,\widehat{g}_0)=\big(b(T_3),g(T_3)\big)$, then
$(\widehat{b}_t,\widehat{g}_t)=(b(T_3+t),g(T_3+t))$ and hence
\begin{equation}\label{equ 5.1}
\widehat{b}_t+\widehat{g}_t\leq \frac{\theta}{2(3+\theta)\lambda}
\end{equation}
for $t\in [0,T_1-T_3]$. Since $\frac{d}{dt}b_t\geq -b_t$, $b_t\geq
b_0e^{-t}$ and hence $b(T_3)>0$. Let $T=T_1-T_3$ and
$\epsilon_1=\min\{\frac{\theta}{5(3+\theta)\lambda},\frac{b(T_3)}{10}\}>0$,
then according to Lemma \ref{Lemma 3.2},
\begin{equation}\label{equ 5.2}
P_{\lambda,n}^{nb(T_3),ng(T_3)}\Big(\sup_{0\leq t\leq
T}\big\|\frac{(B_t,G_t)}{n}-(\widehat{b}_t,\widehat{g}_t)\big\|_1~>\epsilon_1\Big)\leq
\frac{C_1(T_1-T_3,\epsilon_1)}{n}
\end{equation}
for $n\geq N_1(T_1-T_3,\epsilon_1)$, where $C_1$ and $N_1$ are as
defined in Lemma \ref{Lemma 3.2}. By Equations \eqref{equ 5.1} and
\eqref{equ 5.2},
\begin{align}\label{equ 5.3}
&P_{\lambda,n}^{nb(T_3),ng(T_3)}\Bigg(\sup_{0\leq t\leq
T_1-T_3}\frac{\|(B_t,G_t)\|_1}{n}~\leq
\frac{7\theta}{10(3+\theta)\lambda}\\
&\text{~and~}\Big(nb(T_3),ng(T_3)\Big)\succeq
\Big(B(T_1-T_3),G(T_1-T_3)\Big)\Bigg)\geq
1-\frac{C_1(T_1-T_3,\epsilon_1)}{n}\notag
\end{align}
for $n\geq N_1(T_1-T_3,\epsilon_1)$. By Lemma \ref{lemma 2.2} and
Equation \eqref{equ 5.3}, for any $(B,G)\preceq
\big(nb(T_3),ng(T_3)\big)$,
\begin{align}\label{equ 5.4}
&P_{\lambda,n}^{B,G}\Bigg(\sup_{0\leq t\leq
T_1-T_3}\frac{\|(B_t,G_t)\|_1}{n}~\leq
\frac{7\theta}{10(3+\theta)\lambda}\\
&\text{~and~}\Big(nb(T_3),ng(T_3)\Big)\succeq
\Big(B(T_1-T_3),G(T_1-T_3)\Big)\Bigg)\geq
1-\frac{C_1(T_1-T_3,\epsilon_1)}{n}\notag
\end{align}
for $n\geq N_1(T_1-T_3,\epsilon_1)$. By Equation \eqref{equ 5.4} and
utilizing the Markov property for $\sqrt{n}$ times,
\begin{align}\label{equ 5.5}
P_{\lambda,n}^{nb(T_3),ng(T_3)}\Big(\sup_{0\leq t\leq
(T_1-T_3)\sqrt{N}}\frac{\|(B_t,G_t)\|_1}{n}~\leq
\frac{7\theta}{10(\theta+3)\lambda}
\Big)&\geq \Big(1-\frac{C_1(T_1-T_3,\epsilon_1)}{n}\Big)^{\sqrt{n}}\\
&\geq 1-\frac{C_1(T_1-T_3,\epsilon_1)}{\sqrt{n}}\notag
\end{align}
for $n\geq N_1(T_1-T_3,\epsilon_1)$. By Lemma \ref{Lemma 3.2},
\[
P_{\lambda,n}^{n,0}\Big(\|\frac{\big(B(T_1),G(T_1)\big)}{n}-\big(b(T_1),g(T_1)\big)\|_1\leq
\frac{b(T_3)}{10}\Big)\geq 1-\frac{C_1(T_1,\frac{b(T_3)}{10})}{n}
\]
for $n\geq N_1(T_1, \frac{b(T_3)}{10})$ and hence
\begin{equation}\label{equ 5.6}
P_{\lambda,n}^{n,0}\Big(\big(B(T_1),G(T_1)\big)\preceq\big(nb(T_3),ng(T_3)\big)\Big)\geq
1-\frac{C_1(T_1,\frac{b(T_3)}{10})}{n}
\end{equation}
for $n\geq N_1(T_1, \frac{b(T_3)}{10})$. Conditioned on
$\big(B(T_1),G(T_1)\big)\preceq\big(nb(T_3),ng(T_3)\big)$, according
to Lemma \ref{lemma 2.2} and Equation \eqref{equ 5.5},
\[
\sup_{0\leq t\leq
(T_1-T_3)\sqrt{N}}\frac{\|(B_{T_1+t},G_{T_1+t})\|_1}{n}~\leq
\frac{7\theta}{10(3+\theta)\lambda}
\]
with probability at least
$1-\frac{C_1(T_1-T_3,\epsilon_1)}{\sqrt{n}}$ for $n\geq
N_1(T_1-T_3,\epsilon_1)$. Hence by Equation \eqref{equ 5.6},
\[
\sup_{0\leq t\leq
(T_1-T_3)\sqrt{N}}\frac{\|(B_{T_1+t},G_{T_1+t})\|_1}{n}~\leq
\frac{7\theta}{10(3+\theta)\lambda}<\frac{\theta}{(3+\theta)\lambda}
\]
with probability at least
$1-\frac{C_1(T_1-T_3,\epsilon_1)}{\sqrt{n}}-\frac{C_1(T_1,\frac{b(T_3)}{10})}{n}$
for $n\geq \max\{N_1(T_1,
\frac{b(T_3)}{10}),~N_1(T_1-T_3,\epsilon_1)\}$. We choose $N_3$ such
that
\[
\frac{C_1(T_1,\frac{b(T_3)}{10})}{n}\leq \frac{1}{\sqrt{n}}
\]
for $n\geq N_3$, then Lemma \ref{lemma 5.1} holds with
$T_1=T_1(\lambda,\theta)$, $N_2=\max\{N_3, ~N_1(T_1,
\frac{b(T_3)}{10}),~N_1(T_1-T_3,\epsilon_1)\}$, $C_2=T_1-T_3$ and
$C_3=1+C_1(T_1-T_3,\epsilon_1)$.

\qed

To prove Equation \eqref{equ 2.2 main}, we introduce the following
Markov process as an auxiliary model.

Let $\{\widetilde{B}_t\}_{t\geq 0}$ be a continuous-time Markov
process with state space $\{0,1,2,\ldots\}$ and transition rates
function given by
\[
\widetilde{B}\rightarrow
\begin{cases}
\widetilde{B}-1 &\text{~at rate~}\widetilde{B},\\
\widetilde{B}+1 &\text{~at
rate~}\frac{\theta\widetilde{B}}{\theta+3},
\end{cases}
\]
then the following lemma shows that $\widetilde{B}$ dies out in
$O(\log n)$ units of time with high probability conditioned on
$\widetilde{B}_0=n$.
\begin{lemma}\label{lemma 5.2}
\[
P\big(\widetilde{B}_{(1+\frac{\theta}{2})\log
n}=0\big|\widetilde{B}_0=n\big)\geq 1-n^{-\frac{\theta}{2\theta+6}}
\]
for each $n\geq 1$.
\end{lemma}
\proof

Let $h(t)=E\big(\widetilde{B}_t\big|\widetilde{B}_0=n\big)$, then
according to the transition rates function of $\widetilde{B}$,
\[
\frac{d}{dt}h(t)=-h(t)+\frac{\theta
h(t)}{\theta+3}=-\frac{3h(t)}{\theta+3}.
\]
Then, $h(t)=ne^{-\frac{3t}{\theta+3}}$ since $h(0)=n$. By
Chebyshev's inequality,
\[
P\big(\widetilde{B}_{(1+\frac{\theta}{2})\log n}\geq
1\big|\widetilde{B}_0=n\big)\leq h\big((1+\frac{\theta}{2})\log
n\big)=n^{-\frac{\theta}{2\theta+6}}.
\]
Therefore,
\[
P\big(\widetilde{B}_{(1+\frac{\theta}{2})\log
n}=0\big|\widetilde{B}_0=n\big)=1-P\big(\widetilde{B}_{(1+\frac{\theta}{2})\log
n}\geq 1\big|\widetilde{B}_0=n\big)\geq
1-n^{-\frac{\theta}{2\theta+6}}.
\]

\qed

At last we give the proof of Equation \eqref{equ 2.2 main}. From now
on we assume that $\widetilde{B}_0=n$.

\proof[Proof of Equation \eqref{equ 2.2 main}]

Let $T_1$ be as defined in Lemma \ref{lemma 5.1}. On the event
$\|(B_{T_1},G_{T_1})\|_1<\frac{n\theta}{(\theta+3)\lambda}$, we
define
\[
\sigma=\inf\{t\geq
0:~\|(B_{T_1+t},G_{T_1+t})\|_1\geq\frac{n\theta}{(\theta+3)\lambda}\}.
\]
For $t\in [T_1,T_1+\sigma]$, $B\rightarrow B+1$ at rate
\[
\frac{\lambda}{n}BG\leq
\frac{\lambda}{n}B\frac{n\theta}{(\theta+3)\lambda}=\frac{\theta
B}{\theta+3}.
\]
Therefore, conditioned on
$\|(B_{T_1},G_{T_1})\|_1<\frac{n\theta}{(\theta+3)\lambda}$,
\[
B_{T_1+t}\leq \widetilde{B}_{t}
\]
for $t\in [0,\sigma]$ in the sense of coupling. As a result,
\begin{align}\label{equ 5.7}
&P_{\lambda,n}^{n,0}\Bigg(B_{T_1+(1+\frac{\theta}{2})\log
n}=0,\|(B_{T_1},G_{T_1})\|_1<\frac{n\theta}{(\theta+3)\lambda},\sigma>(1+\frac{\theta}{2})\log
n\Bigg)\\
&\geq
P_{\lambda,n}^{n,0}\Bigg(\widetilde{B}_{(1+\frac{\theta}{2})\log
n}=0,\|(B_{T_1},G_{T_1})\|_1<\frac{n\theta}{(\theta+3)\lambda},\sigma>(1+\frac{\theta}{2})\log
n\Bigg)\notag\\
&\geq P\big(\widetilde{B}_{(1+\frac{\theta}{2})\log
n}=0\big)-P_{\lambda,n}^{n,0}\Big(\sup_{0\leq t\leq
(1+\frac{\theta}{2})\log n}\|(B_{T_1+t},G_{T_1+t})\|_1~\geq
\frac{n\theta}{(\theta+3)\lambda}\Big)\notag.
\end{align}
We choose $N_4$ such that $(1+\frac{\theta}{2})\log n\leq
C_2\sqrt{n}$ for each $n\geq N_4$, then by Lemma \ref{lemma 5.1},
\begin{equation}\label{equ 5.8}
P_{\lambda,n}^{n,0}\Big(\sup_{0\leq t\leq (1+\frac{\theta}{2})\log
n}\|(B_{T_1+t},G_{T_1+t})\|_1~\geq
\frac{n\theta}{(\theta+3)\lambda}\Big)\leq \frac{C_3}{\sqrt{n}}
\end{equation}
for $n\geq N_5=\max\{N_4, N_2\}$. By Equations \eqref{equ 5.7},
\eqref{equ 5.8} and Lemma \ref{lemma 5.2},
\begin{equation}\label{equ 5.9}
P_{\lambda,n}^{n,0}\Big(B_{T_1+(1+\frac{\theta}{2})\log
n}=0\Big)\geq 1-\frac{C_3}{\sqrt{n}}-n^{-\frac{\theta}{2\theta+6}}
\end{equation}
for $n\geq N_5$. By Equation \eqref{equ 5.9},
\begin{equation}\label{equ 5.10}
P_{\lambda,n}^{n,0}\Big(\tau\leq T_1+(1+\frac{\theta}{2})\log
n\Big)\geq 1-\frac{C_3}{\sqrt{n}}-n^{-\frac{\theta}{2\theta+6}}
\end{equation}
for $n\geq N_5$. Equation \eqref{equ 2.2 main} follows from Equation
\eqref{equ 5.10} directly since
\[
T_1+(1+\frac{\theta}{2})\log n\leq (1+\theta)\log n
\]
for sufficiently large $n$.

\qed

\quad

\textbf{Acknowledgments.} The author is grateful to Dr. Zhichao
Shan, who suggests us to study the contact process with
semi-infected state. The author is grateful to the financial support
from the National Natural Science Foundation of China with grant
number 11501542 and the financial support from Beijing Jiaotong
University with grant number KSRC16006536.

{}
\end{document}